\documentclass[12pt]{amsart}  

\usepackage[utf8]{inputenc}
\usepackage{amsmath}
\usepackage{amsfonts}
\usepackage{amssymb}
\usepackage{latexsym}
\usepackage{graphicx}
\usepackage{subfigure}
\usepackage{color}
\usepackage{mathtools}
\usepackage[backref=page]{hyperref}
\usepackage{verbatim}
\usepackage[all]{xy}
\usepackage{pdfsync}
\usepackage{youngtab}
\usepackage{ytableau}
\usepackage{young}
\usepackage{tikz}
\usepackage{cancel} 
\usepackage[normalem]{ulem} 
\usepackage{graphicx}
\usepackage{forest}
\usepackage{array}
\usepackage{extarrows}

\theoremstyle{plain}
\newtheorem{theorem}{Theorem}[section]

\newtheorem{corollary}[theorem]{Corollary}

\newtheorem{definition}[theorem]{Definition}
\newtheorem{example}[theorem]{Example}

\newtheorem{remark}[theorem]{Remark}

\numberwithin{equation}{section}
\newcommand{\bN}{\mathbb{N}}

\newcommand{\bQ}{\mathbb{Q}}

\newcommand{\p}[1]{\mathcal #1} 
\newcommand{\smirnov}{\p W}
\newcommand{\N}{{\mathbb N}} %
\DeclareMathOperator{\des}{\mathsf{des}}
\DeclareMathOperator{\asc}{\mathsf{asc}}
\DeclareMathOperator{\rdes}{\mathsf{rdes}}
\DeclareMathOperator{\ldes}{\mathsf{ldes}}
\DeclareMathOperator{\rasc}{\mathsf{rasc}}
\DeclareMathOperator{\lasc}{\mathsf{lasc}}
\newcommand{\alpx}{\mathbf{x}}
\newcommand{\ld}{\lambda}
\newcommand{\la}{\bar{\lambda}}
\newcommand{\ra}{\bar{\rho}}
\newcommand{\rd}{\rho}
\newcommand{\weight}{{\mathsf {wt}}} 
\newcommand{\wordsandtrees}{\p X} 
\newcommand{\smirnovwordgf}{Q} 
\newcommand{\smirnovwords}{\p W} 
\newcommand{\smirnovtrees}{\p T} 
\newcommand{\bleedingtree}{\p B} 
\newcommand{\ppath}{\mathsf{P}} 
\newcommand{\pmax}{\mathsf{M}}
\newcommand{\pmin}{\mathsf{m}} 
\newcommand{\expo}{\mathsf{ex}} 

\definecolor{emphcol}{rgb}{0.0, 0.5, 2.0}
\newcommand{\bemph}[1]{\textcolor{blue}{\emph{#1}}} 
\newcommand{\mbP}{\bN^{*}} 

\newcommand{\Des}[1]{\mathrm{Des}(#1)} 
\newcommand{\Asc}[1]{\mathrm{Asc}(#1)} 

\begin{document}

\title[]{Smirnov Trees}
\author{Matja\v z Konvalinka}
\address{Department of Mathematics, University of Ljubljana \& Institute of Mathematics, Physics and Mechanics, Ljubljana, Slovenia}
\email{\href{mailto:matjaz.konvalinka@fmf.uni-lj.si}{matjaz.konvalinka@fmf.uni-lj.si}}
\author{Vasu Tewari}
\address{Department of Mathematics, University of Pennsylvania, Philadelphia, PA 19104, USA}
\thanks{The first author acknowledges the financial support from the Slovenian Research Agency (research core funding No. P1-0294). The second author was supported by an AMS-Simons travel grant.}
\email{\href{mailto:vvtewari@math.upenn.edu}{vvtewari@math.upenn.edu}}

\begin{abstract}
We introduce a generalization of Smirnov words in  the context of labeled binary trees, which we call Smirnov trees.
We study the generating function for ascent-descent statistics on Smirnov trees and establish that it is $e$-positive, which is akin to the classical case of Smirnov words.
Our proof relies on an intricate weight-preserving bijection.
\end{abstract}

\maketitle
\tableofcontents{}

\section{Introduction}\label{sec:intro}
The study of permutation statistics has been an active area of research since the seminal work of MacMahon \cite{MacMahon}. A particular permutation statistic that plays a prominent role in algebraic combinatorics is the descent statistic.
The associated integer sequence obtained by counting permutations according to their number of descents is the well-known sequence of Eulerian numbers.
Rather than giving an exhaustive list of areas in mathematics where this sequence makes an appearance, we refer the reader to \cite{Petersen} for a beautiful book exposition.
Given that permutations may be considered to be linear trees, it is natural to consider generalizations of classical permutation statistics in the context of labeled plane binary trees.

Gessel \cite{Gessel-Oberwolfach} was the first to study the analogue of the descent statistic for labeled binary trees, and he further pointed out intriguing connections to the enumerative theory of Coxeter arrangements.
There has been a flurry of activity towards understanding these connections better, and the reader is referred to \cite{Bernardi, Corteel-Forge-Ventos, Drake,Gessel-Griffin-Tewari,Tewari,TvW} for more details.
Gessel-Griffin-Tewari \cite{Gessel-Griffin-Tewari} investigated these connections from the perspective of symmetric functions; they attach a multivariate generating function tracking ascent-descents over all labeled binary trees and subsequently prove that this generating function expands positively in terms of ribbon Schur functions. This result is part of the motivation for our article. The rest of it stems from recent work of Shareshian-Wachs \cite{Sh-Wa} on chromatic quasisymmetric functions, which we discuss next.

Stanley \cite{Stanley-chromatic} introduced  chromatic symmetric functions of graphs as a way to generalize  chromatic polynomials of graphs.
In the case where the graph has $n$ nodes and no edges, the chromatic polynomial is the sum of $\alpx_w$ where $w$ ranges over all words of length $n$.
Here, and throughout the rest of this article, set $\alpx_w\coloneqq x_{w_1}\cdots x_{w_n}$ if $w=w_1\dots w_n$, where $\alpx=\{x_1,x_2,\dots\}$ is an alphabet with commuting indeterminates.
It is worth remarking that the refined version tracking the distribution of descents over the set of all words of length $n$  is ribbon Schur-positive.
This motivates studying an analogue of the  descent statistic for general graphs in tandem with the chromatic symmetric function, and this was done by Shareshian and Wachs in \cite{Sh-Wa}, wherein they introduced chromatic quasisymmetric functions.
In the special case where the underlying graph is a path on $n$ nodes, the chromatic quasisymmetric function tracks descents in Smirnov words of length $n$, i.e.~words where two adjacent letters are distinct. 
More importantly, in contrast to the ribbon-positivity in the case where the graph was completely disconnected, the chromatic quasisymmetric function of the path graph is $e$-positive, i.e.~it can be expressed as a non-negative integer linear combination of elementary symmetric functions.
This raises the following natural questions:
\begin{enumerate}
 \item Is there an analogue of Smirnov words and the descent statistic in the context of labeled binary trees?
 \item If yes, is the generating function tracking the distribution of descents $e$-positive?
 \end{enumerate}

In this article, we provide answers to both these questions in the positive by introducing Smirnov trees \textemdash{} labeled rooted plane binary trees with the property that if the parent has the same label as one of its children, then the left child must have a larger label than the right child.
Denote the set of all Smirnov trees by $\smirnovtrees$.
For any labeled binary tree $T$, let $\lasc(T)$, $\ldes(T)$, $\rasc(T)$ and $\rdes(T)$ denote the number of ascents and descents in the labeling to the left and right, such that $\rasc$ and $\lasc$ are determined by weak inequalities, whereas $\rdes$ and $\ldes$ are strict. See Subsection~\ref{subsec:smirnov trees} for details.
We associate a monomial $\alpx_T$ with $T$ as follows. For a node $v\in T$ labeled $i$, let $x_v$ be $x_i$.
Then
\begin{align*}
\alpx_T=\prod_{v\in T}x_v.
\end{align*}
Consider the formal  power series in $\alpx$ with coefficients in $\bQ[\ra,\rd,\la,\ld]$,
\begin{align*}
G:=G(\ra,\rd,\la,\ld)=\sum_{T\in \p T}\ra^{\rasc(T)}\rd^{\rdes(T)}\la^{\lasc(T)}\ld^{\ldes(T)}\alpx_T.
\end{align*}
It is not immediate that $G$ is a symmetric function in $\alpx$ with coefficients in $\mathbb{Q}[\ra,\rd,\la,\ld]$. We establish the preceding fact via the following functional equation satisfied by $G$, which is also our first main result.
\begin{theorem}\label{thm:functional equation smirnov trees}
  Let $\smirnovwords_n$ be the set of Smirnov words of length $n$.
  Then we have
  \[
  G(\ra,\rd,\la,\ld)=\sum_{n\geq 1}\sum_{w\in \smirnov_n}(\ra \la G+\ra+\la)^{\asc(w)}(\rd\ld G+\rd+\ld)^{\des(w)}\alpx_w.
  \]
\end{theorem}
\noindent Theorem~\ref{thm:functional equation smirnov trees} follows from an intricate weight-preserving bijection between the sets  $ \bigcup_{n \geq 1} \smirnovwords_n \times \left( \p T \cup \{D,U\} \right)^{n-1}$ and $ \p T$, and this forms the technical crux of our article. 

Let $e_n\coloneqq e_n(\alpx)$ denote the $n$-th elementary symmetric function, and let $
E(z)\coloneqq\sum_{n\geq 0}e_nz^n$.
It is well known \cite[Theorem C.4]{Sh-Wa} that the generating function tracking ascents and descents over all Smirnov words is $e$-positive.
This fact in conjunction with Theorem~\ref{thm:functional equation smirnov trees} implies the following theorem.
\begin{theorem}\label{thm:e-pos}
$G(\ra,\rd,\la,\ld)$ is $e$-positive.
\end{theorem}
In fact, we provide a combinatorially explicit description for the coefficients in the basis of elementary symmetric functions.
Additionally, we establish the following functional equation satisfied by $G$.
\begin{theorem}\label{thm:another lift of Gessel's}
\begin{align*}
\frac{(1+\ra G)(1+\la G)}{(1+\rd G)(1+\ld G)}=\frac{E(\ra \la G+\ra+\la)}{E(\rd \ld G+\rd+\ld)}.
\end{align*}
\end{theorem}

\medskip

\textbf{Outline of the article:} In Section \ref{sec:backg}, we provide the necessary definitions for Smirnov words and Smirnov trees. In Section \ref{sec:bijection}, we outline how we prove Theorem \ref{thm:functional equation smirnov trees}, how $e$-positivity follows from the functional equation, and how to compute the coefficients in the $e$-expansion.
We give a proof of Theorem~\ref{thm:another lift of Gessel's} as well.
 Finally, in Section \ref{sec:the gory details}, we construct an intermediate bijection. The bijection that proves the main result is then the result of successive applications of this intermediate bijection. We provide many examples and a technical proof.

\section{Background}\label{sec:backg}
We denote the set of positive integers by $\bN$.
For $n\in \bN$, set $[n]\coloneqq \{1,\dots,n\}$.
In particular, $[0]=\emptyset$.
For any undefined terms, we refer the reader to \cite{Stanley-ec2}.

\subsection{Words}\label{subsec:words}
Let $\mbP$ denote the set of all words in the alphabet $\bN$.
Given $w=w_1\dots w_n\in \mbP$, we call $n$ the \bemph{length} of $w$.
An index $i\in [n-1]$ is called a \bemph{descent} of $w$ if $w_i>w_{i+1}$ and an \bemph{ascent} otherwise.
We denote the set of descents (respectively ascents) of $w$ by $\Des{w}$ (respectively $\Asc{w}$) and denote $|\Des{w}|$ (respectively $|\Asc{w}|$) by $\des(w)$ (respectively $\asc(w)$).

A word $w = w_1\dots w_n$ is a \bemph{Smirnov word} if adjacent letters are different, i.e.~if $w_i \neq w_{i+1}$ for $i\in [n-1]$.
The set of all Smirnov words is denoted by $\smirnovwords$ and the set of all Smirnov words of length $n$ by $\smirnovwords_n$.
Consider the multivariate generating function tracking the distribution of ascents and descents over $\smirnovwords_n$
\[
\smirnovwordgf_n(s,t)=\sum_{w\in \smirnov_n}s^{\asc(w)}t^{\des(w)}\alpx_w,
\]
and let $\smirnovwordgf(z;s,t)\coloneqq\sum_{n\geq 0}\smirnovwordgf_n(s,t)z^n$.

\subsection{Smirnov trees}\label{subsec:smirnov trees}
A \bemph{plane binary tree} is a rooted tree in which every node has at most two children, of which at most one is called a \bemph{left child} and at most one is called a \bemph{right child}.
Henceforth, we simply say binary tree instead of plane binary tree.
A \bemph{labeled plane binary tree} (or simply a \bemph{labeled binary tree}) is a binary tree whose nodes have labels drawn from $\bN$.
We assign a weight to labeled binary trees as follows.
An edge between a parent (resp.~a left child) with label $a$ and its right child (resp.~its parent) with label $b$ is weighted $\bar\rho$ (resp.~$\bar\lambda$) if $a \leq b$ and $\rho$ (resp.~$\lambda$) if $a > b$. In other words, if the edge is between a parent and the right child, we use rho (for right), and if it is between a parent and its left child, we use lambda (for left). We add the bar if the nodes form a weak ascent when reading from left to right (either diagonally up or down).
 A node with label $a$ has weight $x_a$.
 The weight of a binary tree $T$, denoted by $\weight(T)$, is the product of the weights of all its edges and nodes.

A labeled binary tree is called \bemph{Smirnov} if the following holds:
\begin{itemize}
\item if the left (resp.~right) child has the same label $i$ as its parent, then the parent must also have a right (resp.~left) child with label $< i$ (resp.~$> i$).
\end{itemize}
In other words, if the parent has the same label as its child, the left child must have a larger label than the right child.
Smirnov trees inherit a weight function from the one  defined for ordinary labeled binary trees.

Denote the set of all Smirnov trees by $\smirnovtrees$, and the set of all Smirnov trees whose root has label $c$ by $\smirnovtrees^c$.
If $T\in \smirnovtrees$, define its \bemph{principal path} $\ppath(T)$ as follows: it starts at the root; if the current node has no right child, stop; if the current node has a left child with the same label (and therefore a right child with a smaller label), move down left; otherwise, move down right.
The last node on the principal path is called the \bemph{principal node} and is denoted $\alpha(T)$; its label is $a(T)$.
The maximal label on $\ppath(T)$ is $\pmax(T)$, and the minimal label is $\pmin(T)$.


\begin{example}
 As an example, let $T$ be the binary tree in Figure~\ref{fig:smirnov-example}.
\begin{figure}[h]
\begin{center}
\includegraphics[]{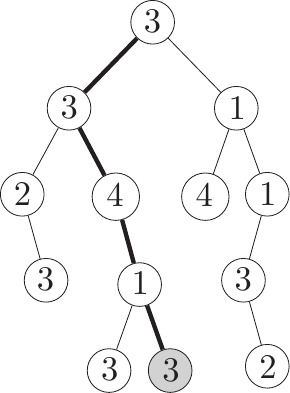}
\caption{A Smirnov tree with principal path in bold.}
\label{fig:smirnov-example}
\end{center}
\end{figure}
There are two nodes with a child with the same label; in both cases, the parent has both children, and the left label is larger than the right label. Therefore $T$ is Smirnov.
In particular, $T \in \smirnovtrees^3$. The edges lying on the principal path $\ppath(T)$ are thickened, and the labels on $\ppath(T)$ are $3,3,4,1,3$. So $a(T) = 3$, $\pmax(T) = 4$, and $\pmin(T) = 1$. The principal node $\alpha(T)$ is gray. Furthermore,
$$\weight(T) =  \bar{\rho}^4 \rho^3 \bar{\lambda}^2 \lambda^3   \, x_1^3x_2^2x_3^6x_4^2.$$
\end{example}

As mentioned in the introduction, we prove Theorem~\ref{thm:functional equation smirnov trees} by establishing a weight-preser\-ving bijection between $\bigcup_{n \geq 1} \smirnovwords_n \times \left( \smirnovtrees \cup \{D,U\} \right)^{n-1}$ and $\smirnovtrees$.
We have already described how to assign weights to elements of the latter, so we proceed to describe weights assigned to elements of the former.
For the sake of brevity, set
\[
\wordsandtrees\coloneqq \bigcup_{n \geq 1} \smirnovwords_n \times \left( \smirnovtrees \cup \{D,U\} \right)^{n-1}.
\]
Given a statement $P$, let $[P]$ equal $1$ if $P$ is true and $0$ otherwise.
If $a,b$ are distinct positive integers and $Y\in \smirnovtrees\cup \{D,U\}$, define
\begin{align*}
\mathsf{f}(a,b,Y)\coloneqq \left\lbrace\begin{array}{ll}
\ra \la\, \weight(T)\, [Y\in \smirnovtrees] +\ra\, [Y=D]+\la \, [Y=U] & \text{ if } a<b,\\
\rd\ld \, \weight(T)\, [Y\in \smirnovtrees] +\rd\, [Y=D]+\ld\, [Y=U]  & \text{ if } a>b.
\end{array}\right.
\end{align*}
This given, consider  $(w,S)\in \wordsandtrees$, and let $w\coloneqq w_1\dots w_n$ and $ S\coloneqq S_1\dots S_{n-1}$.
The weight of $(w,S)$, denoted by $\weight(w, S)$, is defined as follows.
\[
\weight(w, S)\coloneqq \alpx_w\prod_{1\leq i\leq n-1}\mathsf{f}(w_i,w_{i+1}, S_i).
\]
\begin{example}\label{ex:weight wordsandtrees}
Consider the Smirnov word $w = 42534242 \in \p W_8$
and the sequence
$$S =\left(\scalebox{0.3}{\begin{forest}
for tree={circle,draw, l sep=10pt}
[1]
\end{forest}},D,
\scalebox{0.3}{\begin{forest}
for tree={circle,draw, l sep=10pt}
[2]
\end{forest}},D,D,U,
\scalebox{0.3}{\begin{forest}
for tree={circle,draw, l sep=10pt}
[2[2][1]]
\end{forest}}\right) \in \left( \p T \cup \{D,U\}\right)^7
$$
Then we have
\begin{align*}
\mathsf{f}(w_1,w_2,S_1)&=\rd\ld \, x_1\nonumber\\
\mathsf{f}(w_2,w_3,S_2)&=\ra\nonumber\\
\mathsf{f}(w_3,w_4,S_3)&=\rd\ld \, x_2\nonumber\\
\mathsf{f}(w_4,w_5,S_4)&=\ra\nonumber\\
\mathsf{f}(w_5,w_6,S_5)&=\rd\nonumber\\
\mathsf{f}(w_6,w_7,S_6)&=\la\nonumber\\
\mathsf{f}(w_7,w_8,S_7)&=\rd^2\la \ld\, x_1x_2^2.
\end{align*}
It follows that
\[
\weight(w,S)=\ra^2 \rd^5 \la^2 \ld^3 \,x_1^2x_2^6x_3x_4^3x_5.
\]
\end{example}
It is straightforward to see that
\begin{align}\label{eqn:wt generating function smirnov union trees}
\sum_{(w,S)\in \wordsandtrees }\weight(w, S)=\sum_{n\geq 1}\sum_{w\in \smirnov_n}(\ra \la G+\ra+\la)^{\asc(w)}(\rd\ld G+\rd+\ld)^{\des(w)}\alpx_w.
\end{align}
Note that the right hand side in equation~\eqref{eqn:wt generating function smirnov union trees} is exactly the right hand side of the equality in Theorem~\ref{thm:functional equation smirnov trees}.
This explains why it suffices to exhibit a weight-preserving bijection between $\smirnovtrees$ and $\wordsandtrees$ in order to prove Theorem~\ref{thm:functional equation smirnov trees}.

\medskip

\section{Bijection and applications} \label{sec:bijection}

\subsection{Outline of the bijection}
Our weight-preserving bijection between $\smirnovtrees$ and $\wordsandtrees$ uses an intermediate bijection that we call $\Phi$.
Consider a triple $(T,S,b)$, where $S\in\smirnovtrees \cup \{D,U\}$ and $b\in \bN $ is such that $a(T)\neq b$.
Recall that $a(T)$ is the label of the principal node of $T$.
We assign a weight to $(T,S,b)$ as follows:
$$\weight(T,S,b) =x_b\, \weight(T) \mathsf{f}(a(T),b,S).$$
Our map $\Phi$, whose precise description we postpone to Section~\ref{sec:the gory details}, satisfies the following.


\begin{theorem} \label{thm:beast}
 The map $\Phi$ is a well-defined weight-preserving bijection from the set of triples $(T,S,b)$ in $\smirnovtrees \times (\smirnovtrees \cup \{D,U\}) \times \N$ satisfying $a(T) \neq b$ to the set of Smirnov trees with at least 2 nodes, and the principal node of the image has label equal to the third argument. In other words, $\weight(\Phi(T,S,b)) = \weight(T,S,b)$ and $a(\Phi(T,S,b)) = b$. 
\end{theorem}
\noindent The proof of Theorem~\ref{thm:beast}, which is fairly involved, is also presented in Section~\ref{sec:the gory details}.
For the moment, we describe how Theorem~\ref{thm:functional equation smirnov trees} follows from Theorem~\ref{thm:beast}, and subsequently we discuss various applications.
\begin{proof} (of Theorem~\ref{thm:functional equation smirnov trees})
By applying $\Phi$ successively, we construct a weight-preserving bijection
$$\Psi \colon \wordsandtrees \longrightarrow \smirnovtrees,$$
which proves our main theorem.
If $w \in \smirnovwords_1$, let $\Psi(w,\emptyset)$ be the binary tree with a single node labeled $w_1$.
Here $\emptyset$ denotes the empty sequence.
If $n \geq 2$, $w = w_1 \cdots w_n \in \smirnovwords_n$, $S = (S_1,\ldots,S_{n-1}) \in \left( \smirnovtrees \cup \{D,U\} \right)^{n-1}$, define $w' = w_1\cdots w_{n-1}$, $S' = (S_1,\ldots,S_{n-2})$, and
$$\Psi(w,S) = \Phi(\Psi(w',S'),S_{n-1},w_n).$$
We prove by induction on $n$ that for $w \in \smirnovwords_n$ and $S \in \left( \smirnovtrees \cup \{D,U\} \right)^{n-1}$, $\Psi(w, S)$ is well defined and $a(\Psi(w,S)) = w_n$.
For $n = 1$, this is obvious, and if it holds for $n-1$, then $a(\Psi(w',S')) = w_{n-1} \neq w_n$, so $\Phi(\Psi(w',S'),S_{n-1},w_n)$ is well defined, and by Theorem \ref{thm:beast}, $a(\Phi(\Psi(w',S'),S_{n-1},w_n)) = w_n$.

As far as the weight-preserving aspect of $\Psi$ is concerned, note that
\begin{align}\label{eqn:recursive_weight_of_X}
\weight(w,S)
&=
x_{w_n}\cdot
\weight(w',S') \cdot \mathsf{f}(w_{n-1},w_{n},S_{n-1})
\end{align}
Furthermore,
\begin{align}\label{eqn:recursive_weight_of_Psi(X)}
\weight(\Psi(w,S))
&=
\weight(\Phi(\Psi(w',S'),S_{n-1},w_n))\nonumber\\
&=
x_{w_n}\cdot \weight(\Psi(w',S'))\cdot\mathsf{f}(w_{n-1},w_n,S_{n-1}).
\end{align}
In going from the first equality to the second, we have utilized the fact that $\Phi$ is a weight-preserving bijection.
Induction along with a comparison of equations~\eqref{eqn:recursive_weight_of_X} and \eqref{eqn:recursive_weight_of_Psi(X)} implies that $\Psi$ is weight preserving. Bijectivity of $\Psi$ also follows from the bijectivity of $\Phi$.
\end{proof}


\subsection{$e$-positivity and bleeding trees}
To write an explicit expression for $\smirnovwordgf_n(s,t)$, we invoke \cite[Theorem C.4]{Sh-Wa} which states that the generating function for Smirnov words of length $n$ is $e$-positive. More precisely, Shareshian and Wachs prove that
\begin{align} \label{eq:sw}
  \sum_{w \in \smirnovwords_n} t^{\des(w)} \alpx_w = \sum_{m=1}^{\lfloor \frac{n+1}2 \rfloor} \sum_{\stackrel{k_1,\ldots,k_m \geq 2}{\sum k_i = n+1}} e_{(k_1-1,k_2,\ldots,k_m)} t^{m-1} \prod_{i=1}^m [k_i - 1]_t,
\end{align}
where $[a]_t = 1 + t + \cdots + t^{a-1}$ for $a\in \bN$. For example, it is easy to check that
$$\sum_{w \in \smirnovwords_3} t^{\des(w)} \alpx_w = (1+t+t^2) e_3 + t e_{21},$$
which agrees with the formula.

\medskip

It immediately follows from equation \eqref{eq:sw} that
\begin{multline} \label{eq1}
  \smirnovwordgf_n(s,t)= s^{n-1}\sum_{m=1}^{\lfloor \frac{n+1}2 \rfloor} \sum_{\stackrel{k_1,\ldots,k_m \geq 2}{\sum k_i = n+1}} e_{(k_1-1,k_2,\ldots,k_m)} (t/s)^{m-1} \prod_{i=1}^m [k_i - 1]_{t/s}\\
  = \sum_{m=1}^{\lfloor \frac{n+1}2 \rfloor} \sum_{\stackrel{k_1,\ldots,k_m \geq 2}{\sum k_i = n+1}} e_{(k_1-1,k_2,\ldots,k_m)} s^{m-1} t^{m-1} \prod_{i=1}^m (s^{k_i-2} + s^{k_i-3} t + \cdots + t^{k_i-2}) \\
  = \sum_{m=1}^{\lfloor \frac{n+1}2 \rfloor} \sum_{\stackrel{k_1,\ldots,k_m \geq 2}{\sum k_i = n+1}} e_{(k_1-1,k_2,\ldots,k_m)} (s^{k_1-2} + s^{k_1-3} t + \cdots + t^{k_1-2}) \prod_{i=2}^m (s^{k_i-1}t + s^{k_i-2} t^2 + \cdots + s t^{k_i-1})
\end{multline}

Our main result Theorem \ref{thm:functional equation smirnov trees} states that the generating function for Smirnov trees satisfies the functional equation
\begin{equation} \label{eq2}
G= \sum_{n \geq 1} \sum_{w \in \smirnovwords_n} s^{\asc (w)} t^{\des (w)} \alpx_w= Q(1;s,t)-1,
\end{equation}
where
$$s = \bar\rho \, \bar\lambda \, G + \bar\rho + \bar\lambda, \qquad t = \rho \, \lambda \, G + \rho + \lambda.$$

From equations \eqref{eq1} and \eqref{eq2}, it follows that $G$ is $e$-positive. Our goal is to describe $c_\pi=[e_{\pi}]G$, the coefficient of $e_\pi$ in the $e$-expansion of $G$, for an arbitrary partition $\pi$.
We will describe $c_\pi$ in terms of certain trees.
For motivation, let us extract the coefficient at $e_{321}$.

\medskip

Among the terms in \eqref{eq2} (with $s$ and $t$ as above) there are
$$e_{132} \cdot 1 \cdot (s^2 t + s t^2) \cdot (s t) \qquad \mbox{and} \qquad e_{123} \cdot 1 \cdot (s t) \cdot (s^2 t + s t^2) $$
corresponding to compositions $(k_1,k_2,k_3) = (2,3,2)$, $(k_1,k_2,k_3) = (2,2,3)$. For every occurrence of $s$ (resp.~$t$), we can select either $\bar\rho \, \bar\lambda \, G$ (resp.~$\rho \, \lambda \, G$) or $\bar\rho + \bar\lambda$ (resp.~$\rho + \lambda$). Of course, selecting a term with $G$ means that the resulting symmetric function will have degree $> 6$, so the two terms contribute
$$2 ((\bar\rho + \bar\lambda)^2(\rho + \lambda) + (\bar\rho + \bar\lambda)(\rho + \lambda)^2) ((\bar\rho + \bar\lambda)(\rho + \lambda))$$
to $c_{321}$.

\medskip

There are other terms in $c_{321}$, however. For example, among the terms in \eqref{eq1} there is also
$$e_{32} \cdot (s^2 + s t + t^2)(st),$$
corresponding to the composition $(k_1,k_2) = (4,2)$. Since $e_{321} = e_{32} e_1$, the contribution of this term to $c_{321}$ is equal to the coefficient of $e_1$ in $(s^2 + s t + t^2)(st)$. This is equal to $$([e_1] (s^2 + st + t^2))([1](st)) + ([1] (s^2 + st + t^2))([e_1](st)).$$
The function $e_1$ in $s \cdot s + s \cdot t + t \cdot t$ comes from selecting one of the three summands, and then selecting $e_1$ from the first or the second term in the product. Therefore
$$[e_1] (s^2 + st + t^2) = 2 (\bar\rho \, \bar\lambda \, c_1)  (\bar\rho + \bar\lambda) + (\bar\rho \, \bar\lambda \, c_1)  (\rho + \lambda) + (\bar\rho + \bar\lambda)  (\rho \, \lambda \, c_1) + 2 (\rho \, \lambda \, c_1)  (\rho + \lambda).$$
It is clear that $c_1$, the coefficient of $e_1$ in $G$, is $1$. We also have
$$[1] (st) = (\bar \rho + \bar \lambda)(\rho + \lambda),$$
$$[1] (s^2 + st + t^2) = (\bar \rho + \bar \lambda)^2 + (\bar \rho + \bar \lambda)(\rho + \lambda) + (\rho + \lambda)^2,$$
$$[e_1] (st) = \bar\rho \, \bar\lambda \, (\rho + \lambda) + (\bar\rho + \bar\lambda) \, \rho \, \lambda,$$
which gives further terms of $c_{321}$.

\medskip

We can distill the above reasoning into the following definition. Note that it follows from \eqref{eq1} that when extracting a coefficient of $e_\pi$ from some expression, we must pick a part of $\pi$ that will play the role of $k_1$. That is the meaning of red edges below.

\begin{definition}
 A \bemph{bleeding tree} is an unordered rooted tree with the following properties:
 \begin{itemize}
  \item nodes at even depth (including the root) are red and unlabeled, nodes at odd depth are black and labeled with positive integers
  \item every red node has children, and it is connected with a red edge with exactly one of them; all other edges in the tree are black
  \item if a black node is connected to its parent with a black edge, the number of children it has is less than or equal to its label, which cannot be $1$; otherwise, it is strictly smaller than its label (in particular, it has no children if its label is $1$).
 \end{itemize}
 For a partition $\pi$, we denote by $\bleedingtree_\pi$ the (clearly finite) set of bleeding trees whose labels (with repetitions) are precisely the parts of $\pi$.\\
 We give black nodes weights. Assume that a black node has label $r$ and $k$ children.
 \begin{itemize}
  \item If the edge to its parent is black, its weight is
  $$p(r,k)=\sum_{i=1}^{r-1} \sum_{j=0}^i \binom i j \binom{r-i}{r-k-j} (\bar \rho \, \bar \lambda)^{i-j}(\bar \rho + \bar \lambda)^{j} (\rho \, \lambda)^{k-i+j}(\rho + \lambda)^{r-k-j}$$
  \item If the edge to its parent is red, its weight is
  $$\bar p(r,k)=\sum_{i=0}^{r-1} \sum_{j=0}^i \binom i j \binom{r-1-i}{r-1-k-j} (\bar \rho \, \bar \lambda)^{i-j}(\bar \rho + \bar \lambda)^{j} (\rho \, \lambda)^{k-i+j}(\rho + \lambda)^{r-1-k-j}$$
 \end{itemize}
 The weight $\weight(U)$ of a bleeding tree $U$ is defined as the product of the weights of all its black nodes and the number of ways to draw a bleeding tree in the plane so that red edges always come first.
\end{definition}

For example, the bleeding tree with one red node, three black nodes with labels $1$, $2$ and $3$, and a red edge between the red node and the node with label $1$ has weight $2 \bar p(1,0) p(2,0)p(3,0)$, with the factor $2$ coming from drawings

\begin{center}
\scalebox{0.6}{\begin{forest}
for tree={circle,draw, l sep=20pt}
[,fill=red
  [1,edge = {ultra thick,red}
  ]
  [3
  ]
  [2
  ]
]
\end{forest}} \qquad
\scalebox{0.6}{\begin{forest}
for tree={circle,draw, l sep=20pt}
[,fill=red
  [1,edge = {ultra thick,red}
  ]
  [2
  ]
  [3
  ]
]
\end{forest}}.
\end{center}

Note that the requirement that a black node with a black edge to its parent cannot have label $1$ is superfluous, as the weight of such a node would be the empty sum, i.e.~$0$. The requirement about the number of children being smaller than (or equal to) the label is also not needed, as the second binomial coefficient in the sums would be $0$ in that case. In other words, given the weighting of trees, only the first two assumptions are necessary for a tree to be bleeding.

\medskip

\begin{example}
 Figure~\ref{fig:bleeding-1} depicts a bleeding tree in $\bleedingtree_{4332222111}$ with weight
 $$2\bar p(3,2) p(2,1) p(2,1) \bar p(2,0)p(4,2)\bar p(1,0)\bar p(3,0)\bar p(2,0) \bar p(1,0)\bar p(1,0).$$
\begin{figure}[h]
\begin{center}
\includegraphics[]{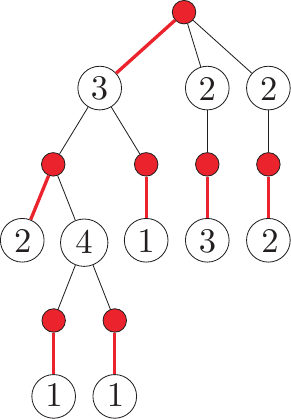}
\caption{A bleeding tree in $\bleedingtree_{4332222111}$.}
\label{fig:bleeding-1}
\end{center}
\end{figure}

All trees in $\bleedingtree_{321}$ are shown in Figure~\ref{fig:bleeding-321}.
\begin{figure}[h]
\begin{center}
\includegraphics[]{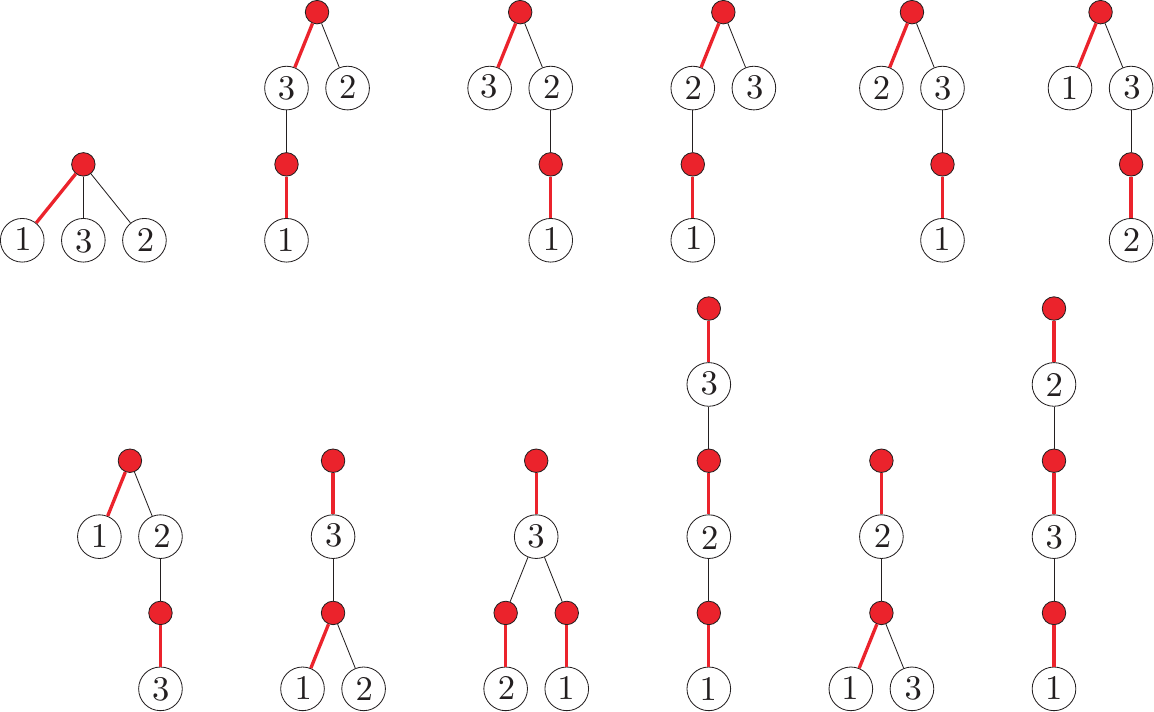}
\caption{All trees in $\bleedingtree_{321}$.}
\label{fig:bleeding-321}
\end{center}
\end{figure}

\end{example}
Given the straightforward yet cumbersome details of the induction, we omit the proof of the following.
\begin{theorem} \label{e-expansion}
 The generating function for Smirnov trees is $e$-positive, with coefficients in $\N[\bar\rho,\rho,\bar\lambda,\lambda]$. More precisely, we can write
 $$G = \sum_{\pi} \left( \sum_{U \in \bleedingtree_\pi} w(U) \right) e_\pi,$$
 where the sum is over all partitions $\pi$, $|\pi| \geq 1$.
\end{theorem}
In fact, the coefficient of $e_{\pi}$ in $G$ is a polynomial with positive integer coefficients in the ``variables'' $\bar \rho \, \bar \lambda$, $\bar \rho + \bar \lambda$, $ \rho \, \lambda$, $ \rho + \lambda$.
Thus, our Theorem~\ref{e-expansion} may be interpreted as an instance of Schur-$\gamma$-nonnegativity defined in \cite{Sh-Wa-2}. See also \cite{Athanasiadis-survey} for more on this theme.
The first few terms of the expansion are
\begin{multline*}
e_1 + (\bar \rho + \rho + \bar\lambda + \lambda) e_2 + (\bar \rho^2 + \bar \rho \rho + \rho^2 + 2 \bar \rho \bar\lambda + \rho \bar\lambda + \bar\lambda^2 + \bar \rho \lambda + 2 \rho \lambda + \bar\lambda \lambda + \lambda^2) e_3 \\+ (\bar \rho \rho + \bar \rho \bar\lambda + \rho \bar\lambda + \bar \rho \lambda + \rho \lambda + \bar\lambda \lambda) e_{21} + \cdots
\end{multline*}
Figure~\ref{fig:bleeding-321} allows us to quickly compute the coefficient at $e_{321}$:
\begin{multline*}
2 \bar p(1,0)p(3,0)p(2,0) + \bar p(3,1) p(2,0) \bar p(1,0) + \bar p (3, 0) p(2, 1) \bar p(1,0) + \bar p(2, 1) p(3, 0) \bar p(1, 0)\\
+ \bar p(2, 0) p(3, 1) \bar p(1, 0) + \bar p(1,0) p(3,1) \bar p(2,0) + \bar p(1,0) p(2,1) \bar p(3,0) + \bar p(3,1) p(2,0) \bar p(1,0) \\
+ 2 \bar p(3,2) \bar p(2,0) \bar p(1,0) + \bar p(3,1) \bar p(2,1) \bar p(1,0) + \bar p(2,1) \bar p(1,0) p(3,0) + \bar p(2,1) \bar p(3,1) \bar p(1,0) \\
= 2 \bar \rho^3 \rho^2 + 2 \bar \rho^2 \rho^3 + 12 \bar \rho^3 \rho \bar \lambda + 18 \bar \rho^2 \rho^2 \bar \lambda + 8 \bar \rho \rho^3 \bar \lambda + 6 \bar \rho^3 \bar \lambda^2 + 28 \bar \rho^2 \rho \bar \lambda^2 + 18 \bar \rho \rho^2 \bar \lambda^2 + 2 \rho^3 \bar \lambda^2 + 6 \bar \rho^2 \bar \lambda^3 \\
+ 12 \bar \rho \rho \bar \lambda^3 + 2 \rho^2 \bar \lambda^3 + 8 \bar \rho^3 \rho \lambda + 18 \bar \rho^2 \rho^2 \lambda + 12 \bar \rho \rho^3 \lambda + 12 \bar \rho^3 \bar \lambda \lambda + 56 \bar \rho^2 \rho \bar \lambda \lambda + 56 \bar \rho \rho^2 \bar \lambda \lambda + 12 \rho^3 \bar \lambda \lambda + 28 \bar \rho^2 \bar \lambda^2 \lambda \\
+ 56 \bar \rho \rho \bar \lambda^2 \lambda + 18 \rho^2 \bar \lambda^2 \lambda + 12 \bar \rho \bar \lambda^3 \lambda + 8 \rho \bar \lambda^3 \lambda + 2 \bar \rho^3 \lambda^2 + 18 \bar \rho^2 \rho \lambda^2 + 28 \bar \rho \rho^2 \lambda^2 + 6 \rho^3 \lambda^2 + 18 \bar \rho^2 \bar \lambda \lambda^2 + 56 \bar \rho \rho \bar \lambda \lambda^2 \\
+ 28 \rho^2 \bar \lambda \lambda^2 + 18 \bar \rho \bar \lambda^2 \lambda^2 + 18 \rho \bar \lambda^2 \lambda^2 + 2 \bar \lambda^3 \lambda^2 + 2 \bar \rho^2 \lambda^3 + 12 \bar \rho \rho \lambda^3 + 6 \rho^2 \lambda^3 + 8 \bar \rho \bar \lambda \lambda^3 + 12 \rho \bar \lambda \lambda^3 + 2 \bar \lambda^2 \lambda^3
\end{multline*}


\subsection{Another functional equation and exponential specialization}\label{subsec:Gessel}
We now give a proof of the functional equation satisfied by $G$ stated in Theorem~\ref{thm:another lift of Gessel's}, and then relate it to earlier work of Gessel.
\begin{proof}(of Theorem~\ref{thm:another lift of Gessel's})
Recall that $\smirnovwordgf(z;s,t)$ tracks the distribution of ascents and descents over the set of all Smirnov words, and we have
\begin{align}\label{eqn:smirnov asc-des}
\smirnovwordgf(z;s,t)&=1+s^{-1}\left( \smirnovwordgf(sz;1,ts^{-1})-1\right).
\end{align}
From \cite[Theorem C.3]{Sh-Wa}, it follows that
\begin{align}\label{eqn:smirnov to elementary}
 \smirnovwordgf(z;1,t)-1=\frac{E(z)-E(zt)}{E(zt)-tE(z)}.
\end{align}
From \eqref{eqn:smirnov to elementary} and \eqref{eqn:smirnov asc-des}, we obtain
\begin{align*}
\smirnovwordgf(1;s,t)-1&=\frac{E(s)-E(t)}{sE(t)-tE(s)}.
\end{align*}
Set $s=\ra \la G+\ra+\la$ and $t=\rd\ld G+\rd+\ld$ henceforth.
By Theorem~\ref{thm:functional equation smirnov trees}, we have
\[
G=\smirnovwordgf(1;s,t)-1,
\]
which in turn implies $(1+\ra G)(1+\la G)=1+Gs=s(\smirnovwordgf(1;s,t)-1)+1$ and $(1+\rd G)(1+\ld G)=1+Gt=t(\smirnovwordgf(1;s,t)-1)+1$.
Thus we obtain
\begin{align}
\frac{(1+\ra G)(1+\la G)}{(1+\rd G)(1+\ld G)}&=\frac{s\frac{E(s)-E(t)}{sE(t)-tE(s)}+1}{t\frac{E(s)-E(t)}{sE(t)-tE(s)}+1}=\frac{E(s)}{E(t)}.
\end{align}
This establishes the claim.
\end{proof}
While the equality in Theorem~\ref{thm:another lift of Gessel's} is less transparent than Theorem~\ref{thm:functional equation smirnov trees}, it immediately allows us to establish a result present in unpublished work of Gessel, and then proved in \cite{Kalikow} and \cite{Drake}.
We call a labeled binary tree on $n$ nodes \bemph{standard} if the labels are all distinct and drawn from $[n]$.
Note that a standard labeled binary trees is necessarily Smirnov.
Gessel considered the following generating function that tracks the distributions of ascents and descents over standard labeled binary trees:
\begin{equation}
B\coloneqq B(\ra,\rd,\la,\ld)=
\sum_{n\geq 1}
\sum_{\substack{T \text{ standard}}}
\ra^{\rasc(T)}\rd^{\rdes(T)}\la^{\lasc(T)}\ld^{\ldes(T)}
\frac{x^n}{n!}.
\end{equation}
Consider the homomorphism $\expo$ from the ring of symmetric functions to $\bQ[[x]]$ defined by setting $\expo(e_n)=\frac{x^n}{n!}$ for  $n\in \bN$.
A key property of $\expo$ is the following. Given a symmetric function $f$, we have  $[x_1\cdots x_n]f=[\frac{x^n}{n!}]\expo(f)$.
It follows that $\expo(G)=B$, and we obtain the following corollary.
\begin{corollary}[]\label{cor:standard_trees}
\begin{equation*}
\frac{(1+\ra B)(1+\la B)}{(1+\rd B)(1+\ld B)}=e^{((\la\ra-\ld\rd)B+\ra+\la-\rd-\ld)x}.
\end{equation*}
\end{corollary}

\section{The gory details}\label{sec:the gory details}
The crucial ingredient of the proof of Theorem~\ref{thm:functional equation smirnov trees} is the map
$$\Phi \colon \left\{ (T,S,b) \in \smirnovtrees \times (\smirnovtrees \cup \{D,U\}) \times \N \colon a(T) \neq b \right\} \longrightarrow \smirnovtrees,$$
which we define as follows. Take a Smirnov tree $T$; $S$, which is either a Smirnov tree, a down step $D$ or an up step $U$; and an integer $b$ that is different from the label of the principal node of $T$.
Write $\alpha = \alpha(T)$, $a = a(T)$, $P = \ppath(T)$, $M = \pmax(T)$, and $m = \pmin(T)$.
By definition, $a \neq b$, and $\alpha$ has no right child.
If $S \in \smirnovtrees$, let $c$ be the label of its root.
If $a < b$ (resp.~$a > b$), let $\delta$ be the last node on $P$ with label $\geq b$ (resp.~$\leq b$), and let $d$ be its label.
If there are no such nodes on $P$, $\delta$ and $d$ are undefined.
For example, if $T$ is the tree from the previous example and $b = 2$, then $d = 1$, and if $b = 5$, $d$ is undefined.
Note that saying, for example, $M < b$ is equivalent to saying that all the labels on the principal path of the tree are smaller than $b$, and $M \geq b$ means that at least one label on the principal path is at least $b$.

\medskip

\begin{figure}[!ht]
\begin{center}
 \includegraphics[scale = 0.7]{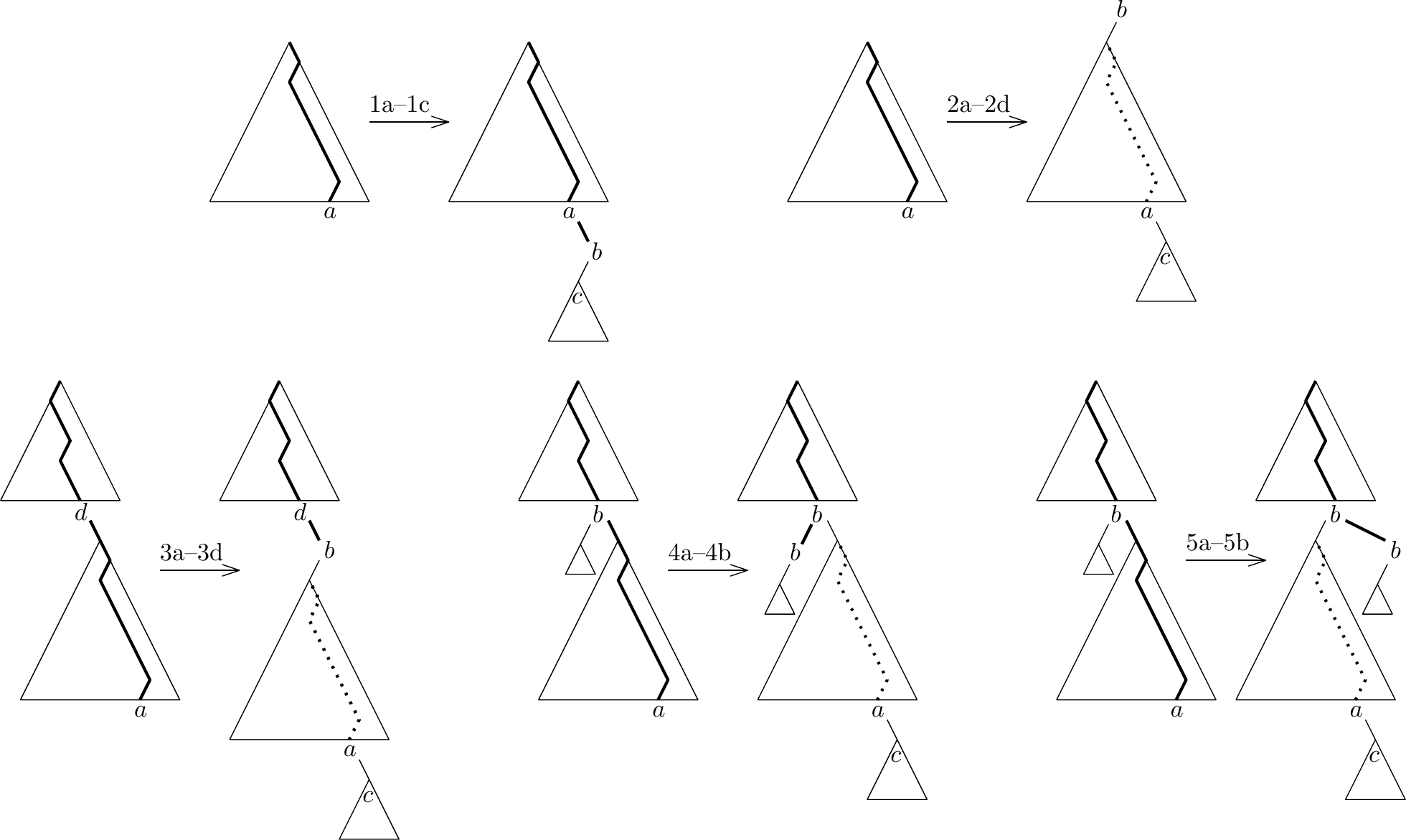}
 \caption{The map $\Phi$.}\label{map}
\end{center}
\end{figure}

The map $\Phi$ is illustrated in Figure \ref{map}, and the formal definition is as follows.

\begin{enumerate}
 \item If
 \begin{itemize}
 \item[(1a)] $S = D$ or
 \item[(1b)] $S \in \smirnovtrees^c$ \& $a,c < b$ or
 \item[(1c)] $S \in \smirnovtrees^c$ \& $a,c>b$,
 \end{itemize}
  then $\Phi(T,S,b)$ is the tree we obtain if we add a right child $\beta$ with label $b$ to $\alpha$.
  If $S = D$, $\beta$ has no children, and if $S \in \smirnovtrees^c$, $\beta$ has $S$ as its left subtree.
 \item If
 \begin{itemize}
 \item[(2a)] $S = U$ \& $a, M < b$ or
 \item[(2b)] $S = U$ \& $a,m > b$ or
 \item[(2c)] $S \in \smirnovtrees^c$ \& $a,M < b \leq c$ or
 \item[(2d)] $S \in \smirnovtrees^c$ \& $a,m > b \geq c$,
 \end{itemize}
  then $\Phi(T,S,b)$ is the tree we obtain if we add a new root $\beta$ with label $b$, and make $T$ its left subtree. If $S \in \smirnovtrees^c$, we make $S$ the right subtree of $\alpha$.
 \item If
 \begin{itemize}
 \item[(3a)] $S = U$ \& $a < b \leq M$ \& $d > b$ or
 \item[(3b)] $S = U$ \& $a > b \geq m$ \& $d < b$ or
 \item[(3c)] $S \in \smirnovtrees^c$ \& $a < b \leq c, M$ \& $d > b$ or
 \item[(3d)] $S \in \smirnovtrees^c$ \& $a > b \geq c, m$ \& $d < b$,
 \end{itemize}
  then $\Phi(T,S,b)$ is the tree we obtain if we replace the right subtree of $\delta$ with a new node $\beta$ with label $b$; $\beta$ has no right child, and its left subtree is the previous right subtree of $\delta$.
  Furthermore, if $S \in \smirnovtrees^c$, we make $S$ the right subtree of $\alpha$.
 \item If
 \begin{itemize}
 \item[(4a)] $S = U$ \& $a < b \leq M$ \& $d = b$ or
 \item[(4b)] $S \in \smirnovtrees^c$ \& $a < b \leq c, M$ \& $d = b$,
 \end{itemize}
 then $\Phi(T,S,b)$ is the tree we obtain if we replace the left subtree of $\delta$ with a new node $\beta$ with label $b$; $\beta$ has no right child, and its left subtree is the previous left subtree of $\delta$.
 Furthermore, if $S \in \smirnovtrees^c$, we make $S$ the right subtree of $\alpha$.
 \item If
 \begin{itemize}
 \item[(5a)] $S = U$ \& $a > b \geq m$ \& $d = b$ or
 \item[(5b)] $S \in \smirnovtrees^c$ \& $a > b \geq c, m$ \& $d = b$,
 \end{itemize} then $\Phi(T,S,b)$ is the tree we obtain if we replace the left subtree of $\delta$ with its right subtree, and its right subtree with a new node $\beta$ with label $b$; $\beta$ has no right child, and its left subtree is the previous left subtree of $\delta$. Furthermore, if $S \in \smirnovtrees^c$, we make $S$ the right subtree of $\alpha$.
\end{enumerate}

\begin{figure}[h]
\begin{center}
\includegraphics[scale=0.9]{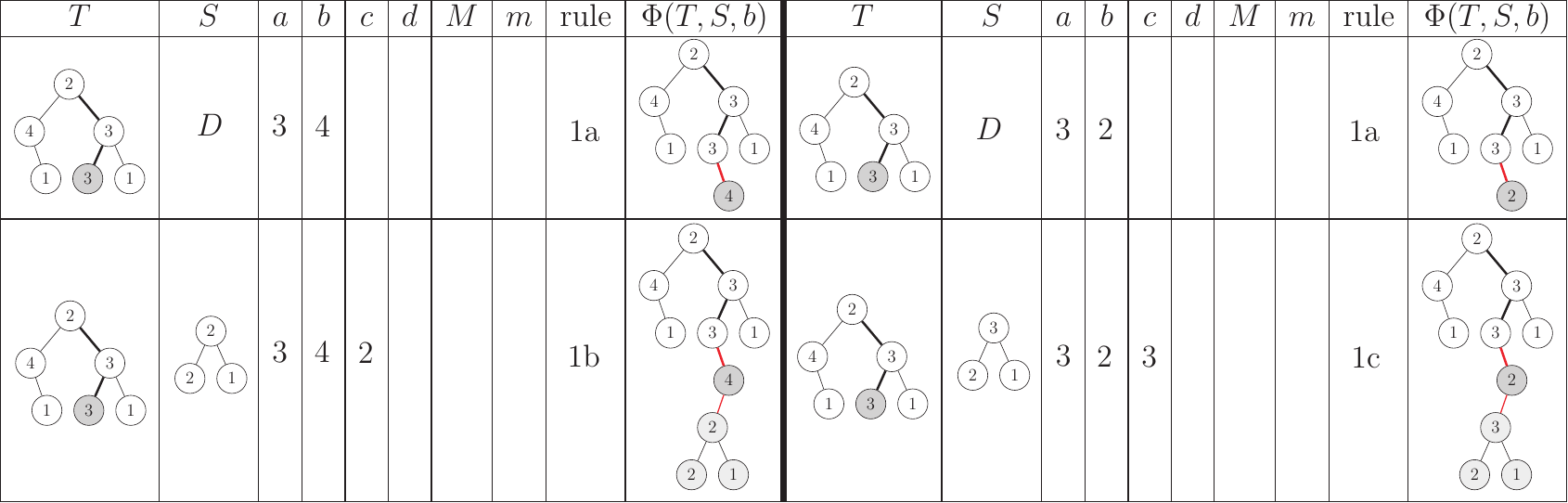}
\caption{An illustration of rule 1.}
\label{fig:demo rule 1}
\end{center}
\end{figure}

\begin{remark}\label{rem:why the bijection is the way it is}
Let us emphasize that while the map is complicated, there is not much freedom if we want to preserve weights and Smirnovness, and if we want the new principal path to end in the node with the newly added label.
\end{remark}
Let us illustrate the simplest rule, rule 1. See Figure~\ref{fig:demo rule 1}.
 We take the same tree $T$ for all four examples, with the principal node $\alpha = \alpha(T)$ being the left child of the right child of the root, and $a = a(T) = 3$ (again, the principal node is gray). By definition, $b$ cannot equal $3$. If $S = D$, we simply add a right child to $\alpha$, with label $b$. The weight of $T$ is multiplied by $x_b$ and, if $a < b$ (resp.~$a > b$), by $\bar\lambda$ (resp.~$\lambda$); the edge that contributes this weight is red. The situation is similar if $S$ is a tree with root label $c$, and $b$ is either larger or smaller than both $a$ and $c$. In that case, add a right child with label $b$ to $\alpha$, and this right child has $S$ as its left subtree (light gray). The weight of $T$ is multiplied by $x_b$, $\weight(S)$, and $\bar\lambda \, \bar\rho$ (resp.~$\lambda \, \rho$) if $a < b$ (resp.~$a > b$); the edges that contribute these weights are red. Note that in all three cases, the principal path is extended by one right step, and $a(\Phi(T,S,b)) = b$.

\begin{figure}[h]
\begin{center}
\includegraphics[scale=0.9]{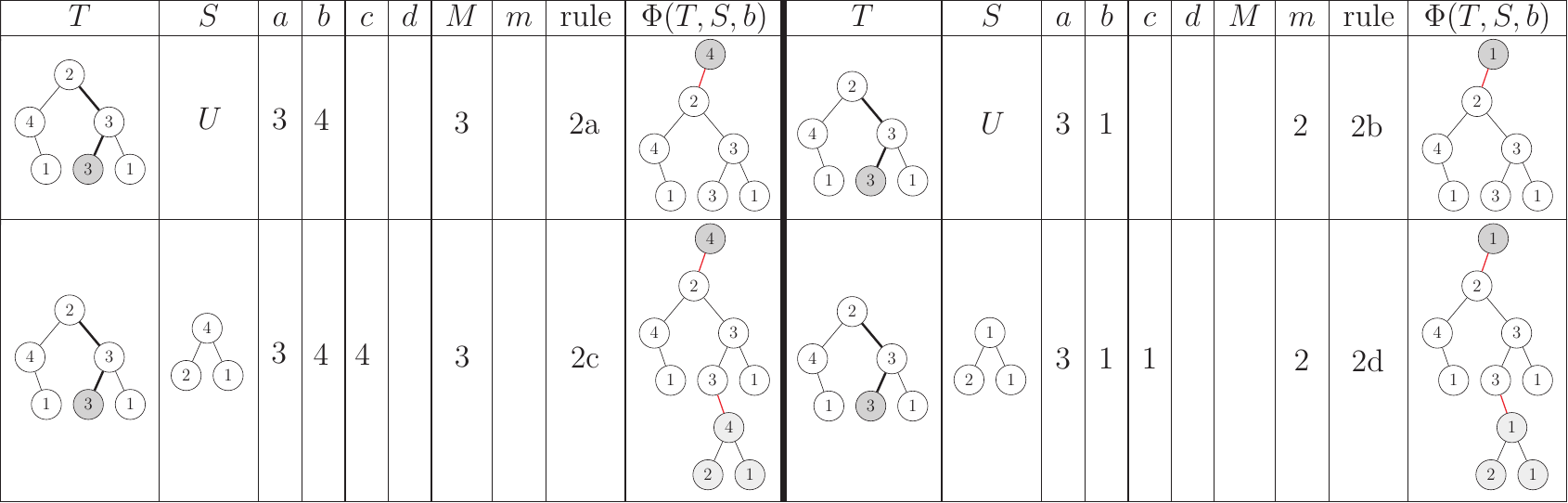}
\caption{An illustration of rule 2.}
\label{fig:demo rule 2}
\end{center}
\end{figure}

Let us present examples of rule 2, which is the only case when we produce a new root. See Figure~\ref{fig:demo rule 2}.
Again, we can take $T$ to be the same tree throughout. If $S = U$, and either $a < b$ and $b$ is larger than all labels on the principal path (i.e.~$M < b$), or $a > b$ and $b$ is smaller than all the labels on the principal path (i.e.~$m > b$), create a new root with label $b$, with no right child, and with the entire $T$ as its left subtree. The subcases 2c and 2d are similar. We have a tree $S$ with root label $c$, and either $a,M < b \leq c$ or $a,m > b \geq c$. In addition to the new root, we attach the tree $S$ as the right subtree of $\alpha$. The new root multiplies the weight of the tree by $x_b$, and if $S$ is a tree, the weight is multiplied by $\weight(S)$. Furthermore, the edge to the new root adds $\bar\lambda$ if $a > b$ and $\lambda$ if $a < b$; if $S$ is a tree, the edge between $\alpha$ and its new right child adds the weight $\bar\rho$ if $a > b$ and $\rho$ if $a < b$.

\medskip

\begin{figure}[h]
\begin{center}
\includegraphics[scale=0.85]{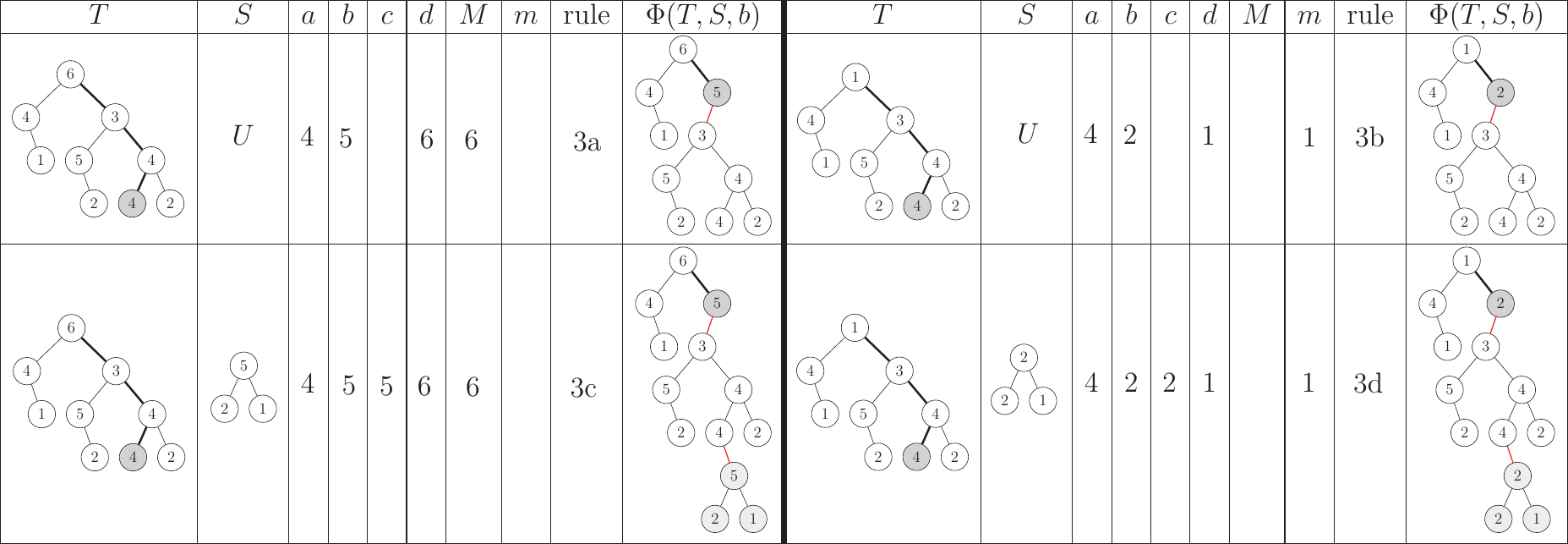}
\caption{An illustration of rule 3.}
\label{fig:demo rule 3}
\end{center}
\end{figure}

 For rule 3, let us take a slightly larger example; note that the $T$'s Figure~\ref{fig:demo rule 3} differ only in the root label. Unlike in rule 2, there is a label on the principal path that is $\geq b$ if $a < b$ and $\leq b$ if $a > b$. Recall that $\delta$ is the last node on the principal path with that property (the root in all our four examples), with label $d$. The crucial assumption now is that $d \neq b$. If $S$ is a tree, its root label must be $\geq b$ if $a < b$ and $\leq b$ if $a > b$ (since otherwise we would use rule 1). Now break up the tree: add a new node $\beta$ as the right child of $\delta$. The new node has no right child, and its left subtree is the old right subtree of $\delta$. Again, if $S$ is a tree, add it as the right subtree of $\alpha$.

Rules 4 and 5 deal with the case when $d = b$. Now $\delta$ becomes a node with the same label as one of its children: left if $a < b$ and right if $a > b$. If $a > b$, the subtrees of $\delta$ are flipped (and we also insert $\beta$ and $S$).

\begin{figure}[h]
\begin{center}
\includegraphics[scale=0.82]{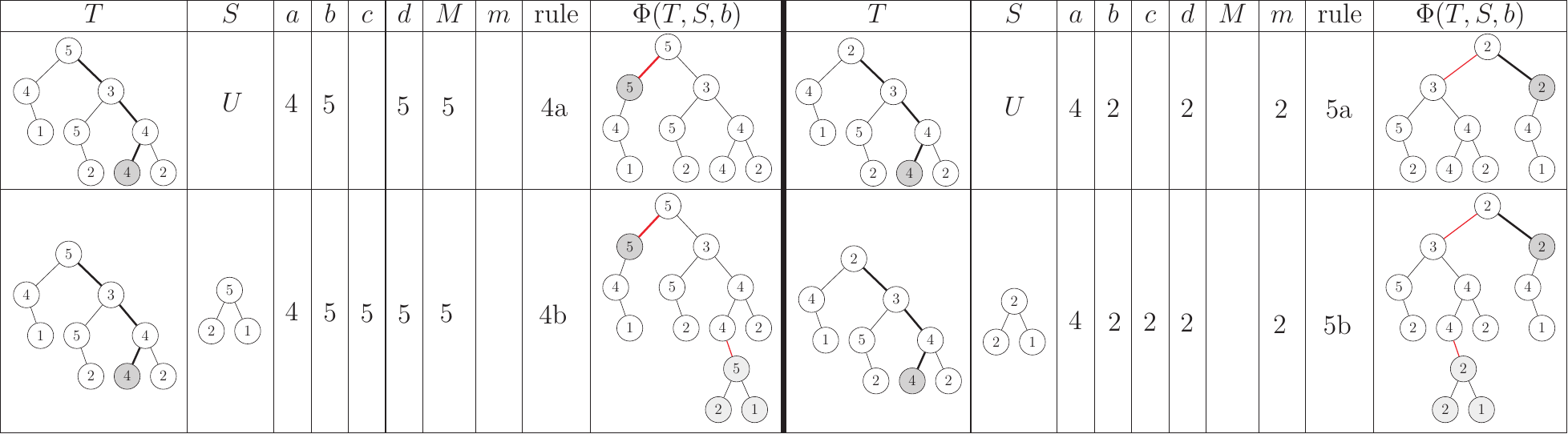}
\caption{An illustration of rules 4 and 5.}
\label{fig:demo rule 45}
\end{center}
\end{figure}

\begin{figure}[h]
\begin{center}
\includegraphics[scale=0.9]{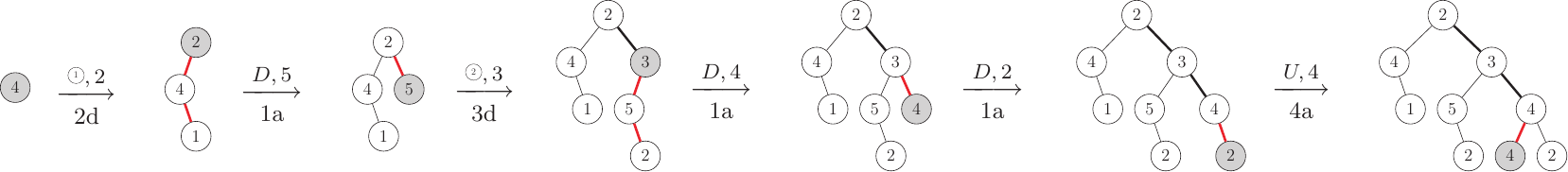}
\caption{$\Phi$ applied to a certain element of $\p X$.}
\label{fig:successive Phi}
\end{center}
\end{figure}

\begin{figure}
\begin{center}
\scalebox{0.55}{\begin{forest}
for tree={circle,draw, l sep=10pt}
[2 [3,
    [5
      [,phantom][2]]
    [4
      [4[,phantom][2[2][1]]][2]]
][{2}[4[,phantom][1]][,phantom]]]
\end{forest}}
\end{center}
\caption{The Smirnov tree corresponding to $w=42534242$ and $\p S$ as in \eqref{eqn:example_S}.}
\label{fig:final_smirnov_tree}
\end{figure}
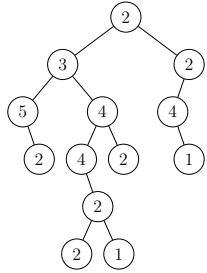

As an example of the use of $\Psi$, consider the Smirnov word
$w = 42534242 \in \p W_8$
and the sequence
\begin{align}\label{eqn:example_S}
\p S = \left(\scalebox{0.3}{\begin{forest}
for tree={circle,draw, l sep=10pt}
[1]
\end{forest}},D,
\scalebox{0.3}{\begin{forest}
for tree={circle,draw, l sep=10pt}
[2]
\end{forest}},D,D,U,
\scalebox{0.3}{\begin{forest}
for tree={circle,draw, l sep=10pt}
[2[2][1]]
\end{forest}}\right) \in \left( \p T \cup \{D,U\}\right)^7.
\end{align}
Then successive applications of $\Phi$ are depicted in Figure~\ref{fig:successive Phi}. The last tree obtained is $T$ from the example for rule (5b), and $S = S_7$ and $b = w_8$ also match that example, so $\Psi(w,\p S)$ is the Smirnov tree in Figure~\ref{fig:final_smirnov_tree}.

\begin{proof}[Proof of Theorem \ref{thm:beast}]
We are given a Smirnov tree $T$; write $\alpha = \alpha(T)$, $a = a(T)$, $M = \pmax(T)$, and $m = \pmin(T)$. We also have $S \in \smirnovtrees \cup \{D,U\}$ and $b \neq a$. We prove that for every appropriate choice of $T,S,b$, exactly one of the rules applies and the result is a Smirnov tree $T'$ with $a(T') = b$ and weight $\weight(T') = \weight(T,S,b)$.\\
The first case is when $S = D$. In that case, rule (1a) applies. Since $a \neq b$ and the principal node $\alpha$ has no right child by definition, we can give it a right child $\beta$ with label $b$. This creates a new right edge in the tree; if $a < b$, then its weight is $\bar\rho$, and if $a > b$, then its weight is $\rho$. In both cases, the weight of the new tree, which is clearly Smirnov, equals the weight of $(T,S,b)$, which is $\bar \rho \,\weight(T) x_b$ or $\rho \, \weight(T) x_b$, depending on whether $a < b$ or $a > b$. It is also clear that the principal path is the same as before with an added right step at the end. So the new principal node has label $b$, as claimed. \emph{Note that in case (1a), the new principal node has no children, and its left parent has a different label. (Statements in italics will be crucial in the construction of the inverse of $\Phi$.)}\\
Now assume that $S = U$ and $a < b$.
The principal path contains a node with label $< b$, and either $M$, the maximum of all labels on the principal path, is $< b$, or it is $\geq b$.
The first case gives rule (2a).
Indeed, since $M < b$, the root label is also $< b$.
When we add a new root $\beta$ with label $b$ and make $T$ its left subtree, we create a new edge with weight $\bar\lambda$, so the weight changes as claimed.
The new tree is also obviously Smirnov.
The principal path contains only the root, and the new principal node again has label $b$. \emph{Note that in case (2a), the root of the new tree is the principal node with a left child, and all the labels on the principal path of its left subtree are smaller than the root label.}\\
The next case is when $S = U$, $a < b$, and $M \geq b$.
On the principal path, there is at least one node with label at least $b$.
Let $\delta$ be the last node on the principal path with that property, and let $d \geq b$ be its label.
The children of $\delta$ must have a label different from $d$ (otherwise, the node after $\delta$ on the principal path is its child with the same label, and therefore $\delta$ is not the last node on the principal path with label $\geq b$).
That means that the principal path goes right from $\delta$.
Again, we have two options, $d > b$ and $d = b$.
The first case gives rule (3a).
Since $d \neq b$, we can replace the right child of $\delta$ with a new node $\beta$ with label $b$, and make the old right subtree of $\delta$ the left subtree of the new node (which has no right child).
Since $d > b$, the edge between $\delta$ and its new right child has label $\rho$; since all the labels on the principal path of the right subtree of $\delta$ in $T$ are $< b$ by the selection of $\delta$, the old edge between $\delta$ and its right child had weight $\rho$, and the new edge between $\beta$ and its left child has label $\bar\lambda$.
In other words, the weight of the new tree is the previous weight times $\bar\lambda \, x_b$, as required.
Also, the new principal path again ends at $\beta$, which has label $b$. \emph{Note that in case (3a), the new principal node is a right child whose parent has a larger label, and all the labels on the principal path of its (nonempty) left subtree are smaller.}\\
The other option is that $S = U$, $a < b$, $M \geq b$, and $d = b$.
This gives rule (4a): we insert a new node $\beta$ with label $b$ on the edge between $\delta$ and its left child (if $\delta$ has no left child, we just add a left child with label $b$).
Since the right child of $\delta$ has label $< b$, and the left child of $\delta$ has label $\neq d$, the new tree is Smirnov.
We added an edge between $\beta$ and $\delta$, which has, since both have labels $b$, weight $\bar\lambda$.
The new principal path now goes to the left after $\delta$, and stops there.
The label of the new principal node is $b$. \emph{Note that in case (4a), the new principal node is a left child, and its label is larger than all the labels on the principal path of its parent's (nonempty) right subtree.}\\
This completes the analysis of the case $S = U$, $a < b$.
The analysis for $S = U$, $a > b$ is similar (the major difference being in the case $d = b$) and we will provide fewer details.
The principal path contains a node with label $> b$, and $m$, the minimum of all labels on the principal path, can be $> b$ or $\leq b$.
The first case gives rule (2b): add a new root, and make $T$ its left subtree.
The new edge has label $\lambda$, so the weight changes as required. \emph{Note that in case (2b), the root is the principal node, and all the labels on the principal path of its (nonempty) left subtree are larger than the root label.}\\
If $S = U$, $a > b$, and $m \leq b$.
On the principal path, there is at least one node with label at most $b$.
Let $\delta$ be the last node on the principal path with that property, and let $d \leq b$ be its label.
Again, the children of $\delta$ must have a label different from $d$, and that means that the principal path goes right from $\delta$.
Again, we have two options, $d < b$ and $d = b$. The first case gives rule (3b).
Replace the right child of $\delta$ with a new node $\beta$ with label $b$, and make the old right subtree of $\delta$ the left subtree of the new node (which has no right child).
The new tree is Smirnov and its weight is the previous weight times $\lambda \, x_b$.
The new principal path again ends at $\beta$, which has label $b$. \emph{Note that in case (3b), the new principal node is a right child whose parent has a smaller label, and all the labels on the principal path of its (nonempty) left subtree are larger.}\\
On the other hand, if $S = U$, $a > b$, $m \leq b$, and $d = b$, we have rule (5a).
We insert a new node $\beta$ with label $b$ as a right child of $\delta$ (not left child as in case (4a)!), make its left subtree the old \emph{left} subtree of $\delta$, and make the old right subtree of $\delta$ its left subtree.
Since the root of the old right subtree of $\delta$ has label $> b$, this is still a Smirnov tree.
The difference in weights comes from the edge between $\delta$ and its new left child, which has weight $\lambda$, and, as always, from the new node $\beta$.
The new principal path ends at $\beta$. \emph{Note that in case (5a), the new principal node is a right child, and its label is the same as its parent's label and smaller than all the labels on the principal path of its parent's (nonempty) left subtree. The new principal node might have no children (if $\delta$ had no left child in $T$).}\\
This leaves us with the case when $S$ is a tree with, say, root label $c$.
The analysis is similar and we will be brief.\\
The first case is when $c$ is ``on the same side'' of $b$ as $a$, i.e.~if $a,c < b$ (use (1b)) or $a,c > b$ (use (1c)).
Apart from adding the new node $\beta$ as the right child of $\alpha$, we also add the tree $S$ as the left subtree of $\beta$.
If $a < b$, the new edges have weights $\bar\rho$ and $\bar\lambda$, and if $a > b$, the weights are $\rho$ and $\lambda$.
Of course, we also added the tree $S$, which contributes $\weight(S)$ to the weight of the new tree.
Thus in both cases, the weight of the new tree equals $\weight(T,S,b)$, and the new principal node is $\beta$. \emph{Note that in case (1b), the new principal node is a right child and has a left child, and its label is larger than the labels of its parent and its left child. In case 1c, the new principal node is a right child and has a left child, and its label is smaller than the labels of its parent and its left child.}\\
Otherwise, $c$ is ``on the other side'' of $b$ as $a$: we have either $a < b \leq c$ or $a > b \geq c$.
Let us assume that $a < b \leq c$.
Again, we look at the labels on the principal path of $T$. If they are all smaller than $b$, we have case (2c).
Apart from adding the new root and making $T$ its left subtree, we also add $S$ as the right subtree of $\alpha$. The new edges have weights $\bar\lambda$ (since the principal label is $< b$) and $\bar \rho$ (since $a < c$), which ensures that the weight of the new tree is $\weight(T,S,b)$. \emph{Note that in case (2c), the root is the principal node, its label is larger than the label of its left child, but there is a node on the principal path of its left subtree (the root of $S$!) that has a label at least as large.}\\
We can also have nodes on the principal path that have labels $\geq b$.
Again, let $\delta$ be the last such node, and let $d \geq b$ be its label (different from the nodes of $\delta$'s children).
If $d > b$, we have rule (3c): do the same as in case (3a), but also add the tree $S$ as the right subtree of $\alpha$. \emph{In case (3c), the new principal node is a right child whose parent has a larger label, the left child has a smaller label, and there is a node on the principal path of its left subtree with a label that is at least as large.}\\
If, however, $d = b$, we are in case (4b), which is similar to (4a). \emph{In case (4b), the new principal node is a left child and its label is larger than the label of its parent's right child, but there is a label on the principal path of its parent's right subtree with a label that is at least as large.}\\
This leaves us with the analysis of the cases when $S \in \smirnovtrees^c$ and $a > b \geq c$.
If $m > b$, we have case (2d), which has the same construction as (2c), but new edge weights $\lambda$ and $\rho$ instead of $\bar\lambda$ and $\bar\rho$.
Of course, this is consistent with the requirement that $a > b$, $S \in \smirnovtrees$ adds weight $\lambda \, \rho \, w(S) x_b$. \emph{In case (2d), the root is the principal node, its label is smaller than the label of its left child, but there is a node on the principal path of its left subtree that has a label at most as large.}\\
If $m \leq b$ and $d < b$, we have case (3d).
Again, do the same thing as in case (3c). \emph{Note that in case (3d), the new principal node is a right child whose parent has a smaller label, the left child has a larger label, and there is a node on the principal path of its left subtree with a label that is at most as large.}\\
Finally, consider the case $m \leq b$, $d = b$, which is rule (5b).
The situation is similar as in (5a), except that we add the tree $S$ as the right child of $\alpha$.
\emph{In case (5b), the new principal node is a right child, its label is the same as its parent's label and smaller than the label of its parent's left child, but there is a label on the principal path of its parent's left subtree with a label that is at most as large. The new principal node might not have children.}\\
This concludes the proof that $\Phi$ is a well-defined weight-preserving map.
The proof is finished if we construct the inverse.
The detailed analysis of the previous paragraphs makes that easy.\\
We are given a tree $T'$ with at least two nodes, and we want to see that there is a unique triple $(T,S,b)$ with $\Phi(T,S,b) = T'$.
Clearly, $b$ must be $a(T')$; denote $\alpha(T')$ by $\beta$.\\
As the principal node, $\beta$ cannot have a right child.
Let us first assume that it has no parent, i.e.~that it is the root of $T'$.
Since $T'$ is not a single node, it must have a nonempty left subtree $T$.
We obtain $T'$ by applying one of the two rules to $T$ or to $T$ with a subtree removed.
More precisely, if all the labels on the principal path of the left subtree of the root are smaller (resp.~larger) than $b$, then $T' = \Phi(T,U,b)$ is the result of rule (2a) (resp.~(2b)).
If, however, the left child of the root has a label that is smaller (resp.~larger) than $b$, but there is a node on the principal path of $T$ whose label is $\geq b$ (resp.~$\leq b$), denote the first such node $\gamma$.
Then we obtain $T$ by applying rule (2c) (resp.~(2d)) to $(T\setminus S,S,b)$, where $S$ is the subtree of $T$ with root $\gamma$.\\
Now assume that $\beta$ is not a root.
The first case is that it is a left child.
By construction of the principal path and the fact that $T'$ is Smirnov, that must mean that its parent's label is also $b$, and that the parent has a right child with a smaller label.
If all the labels on the principal path of $\beta$'s parent's right subtree are smaller than $b$, we obtain $T'$ by applying rule (4a) to $(T,U,b)$, where $T$ is obtained from $T$ by deleting the node $\beta$ (if $\beta$ has a left subtree, simply make this subtree $\beta$'s parent's left subtree).
If, on the other hand, $\gamma$ is the first node on the principal path of $\beta$'s parent's right subtree with label $\geq b$ and $S$ is the subtree of $T'$ with root $\gamma$, then $T' = \Phi(T \setminus S,S,b)$.\\
The last (and most common, as the principal path prefers the right direction) case is when $\beta$ is the right child. The first case is when its parent's label is also $b$.
That must necessarily mean that $\beta$'s parent has a left child with label $> b$.
If all the labels on the principal path of $\beta$'s parent's left subtree are $> b$, we obtain $T'$ as an application of rule (5a) to $(T,U,b)$, where $T$ is the tree we obtain if we remove $\beta$ and swap the resulting left and right trees of $\beta$'s parent. Otherwise, $T'$ results by using rule (5b).\\
Now we can assume that $\beta$ is the right child, and its parent has a different label.
If $\beta$ has no children, we can produce it by using rule (1a) on $(T' \setminus \{\beta\},D,b)$. If $\beta$ has a left child, it will necessarily have a different label. If the labels of $\beta$'s parent and child are both smaller than $b$, $T'$ is obtained by rule (1b), and if they are both larger than $b$, by rule (1c).\\
That leaves us with the case when one of the labels is larger and the other one is smaller.
If the parent's label is larger (resp.~smaller) and all the labels on the principal path of $\beta$'s left subtree are smaller (resp.~larger), we applied rule (3a) (resp.~(3b)) to get $T'$.
If, however, the parent's label is larger (resp.~smaller), the left child's label is smaller (resp.~larger), but there is a node on the principal path of $\beta$'s child whose label is at least $b$ (resp.~at most $b$), we applied rule (3c) (resp.~(3d)).
This completes the construction of the inverse, and the proof of the theorem.
\end{proof}

A suggested exercise for the reader is to take the tree drawn just before the proof and reproduce the word $w$ and the sequence $\p S$.

\section{Final remarks}
\begin{enumerate}
\item As mentioned in the introduction, Smirnov words can be interpreted as proper colorings of path graphs. Smirnov trees, on the other hand, are a labeled tree-analogue of Smirnov words.
This raises the natural question of constructing labeled (binary) tree analogues for graphs other than path graphs, and see if one can define ascent-descent statistics on these tree analogues that relate to the Shareshian-Wachs chro\-mat\-ic quasisymmetric function.
Another potential avenue is to consider the directed or cyclic analogue of the aforementioned question given the recent work of Ellzey-Wachs \cite{Ell-Wa} and Panova-Alexandersson \cite{Pan-Al}.
\item  As established in Subsection~\ref{subsec:Gessel}, the exponential specialization of $G$ tracks ascents-descents in standard labeled binary trees.
In particular, this implies the  equalto preserveities
\begin{align*}
[x_1\dots x_n]G(1,1,1,1)&=n!\times \mathrm{Cat}_{n},\\
[x_1\dots x_n]G(1,1,1,0)&=(n+1)^{n-1},
\end{align*}
where $\mathrm{Cat}_n$ denotes the $n$th Catalan number.
It is unclear to us how one could derive these equalities starting from our $e$-positive expansion for $G$ in terms of weights of bleeding trees.
\item The $h$-positivity of $\omega G$ (here $\omega$ is the standard involution on the algebra of symmetric functions) raises the question of constructing a natural permutation rep\-re\-sen\-ta\-tion that realizes $\omega G$ as its Frobenius characteristic.
Stanley \cite[Proposition 7.7]{Stanley}, using a recurrence due to Procesi \cite{Pr}, established that $\omega G(\ra,\rd,0,0)$ can be realized as the generating function of the Frobenius characteristic of the representation of $S_n$ on the cohomology of the toric variety associated with the Coxeter complex of type $A_{n-1}$.
Stembridge \cite{Stembridge} also constructed a symmetric group representation which realizes $\omega G(\ra,\rd,0,0)$ as the Frobenius char\-ac\-ter\-is\-tic via a bijection from permutations to what he calls codes.
It would be interesting to generalize these results to Smirnov trees.

The following are the characters that such permutation representations would have. For example, there are 288 standard trees on five nodes that have weight $\ra^2 \la \ld$ (i.e., that have two right ascents, one left ascent and one left descent). According to line 7 of the last table, 94 of them should be fixed under the action of a transposition, 8 under the action of a $4$-cycle etc.

$$
\begin{array}{cccc}
 \text{} & 111 & 21 & 3 \\
  \ra^2, \rd^2, \la^2, \ld^2 & 1 & 1 & 1 \\
 \ra \rd, \ra \ld, \rd \la, \la \ld & 4 & 2 & 1 \\
 \ra \la, \rd \ld & 5 & 3 & 2 \\
\end{array} \qquad
\begin{array}{cccccc}
 \text{} & 1111 & 211 & 22 & 31 & 4 \\
  \ra^3, \rd^3, \la^3, \ld^3 & 1 & 1 & 1 & 1 & 1 \\
 \ra^2 \rd, \ra^2 \ld, \ra \rd^2, \ra \ld^2, \rd^2 \la, \rd \la^2, \la^2 \ld,  \la \ld^2 & 11 & 5 & 3 & 2 & 1 \\
 \ra^2 \la, \ra \la^2, \rd^2 \ld, \rd \ld^2 & 17 & 9 & 5 & 5 & 3 \\
 \ra \rd \la, \ra \rd \ld, \ra \la \ld, \rd \la \ld & 44 & 16 & 8 & 5 & 2 \\
\end{array}$$

$$\begin{array}{cccccccc}
 \text{} & 11111 & 2111 & 221 & 311 & 32 & 41 & 5 \\
  \ra^4, \rd^4, \la^4, \ld^4 & 1 & 1 & 1 & 1 & 1 & 1 & 1 \\
 \ra^3 \rd, \ra^3 \ld, \ra \rd^3, \ra \ld^3, \rd^3 \la, \rd \la^3, \la^3 \ld, \la \ld^3 & 26 & 12 & 6 & 5 & 3 & 2 & 1 \\
 \ra^3 \la, \ra \la^3, \rd^3 \ld, \rd \ld^3 & 49 & 25 & 13 & 13 & 7 & 7 & 4   \\
 \ra^2 \rd^2, \ra^2 \ld^2, \rd^2 \la^2, \la^2 \ld^2 & 66 & 22 & 10 & 6 & 4 &   2 & 1 \\
 \ra^2 \la^2, \rd^2 \ld^2 & 146 & 60 & 26 & 26 & 12 & 12 & 6 \\
 \ra^2 \rd \ld, \ra \rd^2 \la, \ra \la \ld^2, \rd \la^2 \ld & 237 & 73 & 29 & 18 & 10 & 5 & 2 \\
 \ra^2 \rd \la, \ra^2 \la \ld, \ra \rd^2 \ld, \ra \rd \la^2, \ra \rd \ld^2, \ra \la^2 \ld, \rd^2 \la \ld, \rd \la \ld^2 & 288 & 94 & 36 & 27 & 13 & 8 & 3 \\
 \ra \rd \la \ld & 824 & 228 & 80 & 50 & 24 & 12 & 4 \\
\end{array}
$$
\end{enumerate}

\section*{Acknowledgements}

The authors would like to thank Guillaume Chapuy, Ira Gessel, Ian Goulden, Sean Griffin, David Jackson, Andrea Sportiello, Stephan Wagner, and Philip Zhang for interesting conversations.


\begin{thebibliography}{00}




\bibitem{Athanasiadis-survey}
{\sc C.~Athanasiadis},
{\em{Gamma-positivity in combinatorics and geometry}},
S\'em. Lothar. Combin. 77 ([2016-2018]), Art. B77i, 64 pp.


\bibitem{Bernardi}
{\sc O.~Bernardi},
{\em{Deformations of the braid arrangement and trees}},
Adv. Math. 335 (2018), 466--518.


\bibitem{Corteel-Forge-Ventos}
{\sc S.~Corteel, D.~Forge and V.~Ventos},
{\em{Bijections between affine hyperplane arrangements and valued graphs}},
European J. Combin. 50 (2015), 30--37.


\bibitem{Drake}
{\sc B.~Drake},
{\em{An inversion theorem for labeled trees and some limits of areas under lattice paths}},
Phd. Dissertation, Brandeis University (2008).

\bibitem{Ell-Wa}
{\sc B.~Ellzey and M.~Wachs},
{\em{On enumerators of Smirnov words by descents and cyclic descents}},
arxiv: \url{https://arxiv.org/abs/1901.01591}.

\bibitem{Gessel-Oberwolfach}
{\sc I.~Gessel},
{\em{Oberwolfach Reports (Enumerative Combinatorics)}},
\url{https://www.mfo.de/document/1410/OWR_2014_12.pdf}, 2014, Page 709.

\bibitem{Gessel-Griffin-Tewari}
{\sc I.~Gessel, S.~Griffin and V.~Tewari},
{\em{Labeled plane binary trees and Schur-positivity}},
arxiv: \url{https://arxiv.org/abs/1706.03055}.


\bibitem{Kalikow}
  {\sc L.~Kalikow},
  {\em {Symmetries in trees and parking functions}},
  Discrete Math. 256 (2002), 719--741.


\bibitem{MacMahon}
{\sc P.A.~MacMahon},
{\em{Combinatory Analysis}}, 2 volumes, Cambridge University Press, London, 1915-1916.


\bibitem{Pan-Al}
{\sc G.~Panova and P.~Alexandersson},
{\em{LLT polynomials, chromatic quasisymmetric functions and graphs with cycles}},
Disc. Math. 341 (2018), 3453--3482.

\bibitem{Petersen}
{\sc T.~Kyle Petersen},
{\em {Eulerian numbers}},
Birkha\"user Advanced Texts: Basler Lehrb\"ucher, Birkha\"user/Springer, New York, 2015.


\bibitem{Pr}
{\sc C.~Procesi},
{\em{The toric variety associated to Weyl chambers}},
Mots, Lang. Raison. Calc., Herm\'es, Paris, 1990, 153--161.

\bibitem{Sh-Wa}
{\sc J.~Shareshian and M.~Wachs},
{\em{Chromatic quasisymmetric functions}},
Adv. Math. 295 (2016), 497--551.

\bibitem{Sh-Wa-2}
{\sc J.~Shareshian and M.~Wachs},
{\em{Gamma-positivity of variations of Eulerian polynomials}},
arxiv: \url{https://arxiv.org/abs/1702.06666}

\bibitem{Stanley}
{\sc R.~Stanley},
{\em{Log-concave and unimodal sequences in algebra, combinatorics, and geometry. }}
Ann. New York Acad. Sci. 576 (1989), 500--535.

\bibitem{Stanley-chromatic}
{\sc R.~Stanley},
{\em{A symmetric function generalization of the chromatic polynomial of a graph}},
Adv. Math. 111 (1995), 166--1994.

\bibitem{Stanley-ec2}
{\sc R.~Stanley},
 {\em {Enumerative Combinatorics} vol.~2},
 Cambridge University Press, 1999.

\bibitem{Stembridge}
{\sc J.~Stembridge},
{\em{Eulerian numbers, tableaux, and the Betti numbers of a toric variety}},
Disc. Math. 99 (1992), 307--320.

\bibitem{Tewari}
{\sc V.~Tewari},
{\em{Gessel polynomials, rooks, and extended Linial arrangements}},
J. Combin. Theory Ser. A 163 (2019), 98--117.

\bibitem{TvW}
{\sc V.~Tewari and S.~van Willigenburg},
{\em{Permuted composition tableaux, 0-Hecke algebra and labeled binary trees}},
J. Combin. Theory Ser. A 161 (2019), 420--452.

\end{thebibliography}
\end{document}